\documentclass[12pt,oneside]{amsart}
\usepackage{amsthm,amsfonts,amssymb}
 \makeatletter
\@addtoreset{equation}{section}\makeatother

\newtheorem{theorem}[subsection]{Theorem}
\newtheorem{proposition}[subsection]{Proposition}
\newtheorem{lemma}[subsection]{Lemma}
\newtheorem{claim}[subsubsection]{Claim}
\newtheorem*{sublemma*}{Sublemma}
\newtheorem{corollary}[subsubsection]{Corollary}
\newtheorem{conjecture}[subsection]{Conjecture}
\theoremstyle{definition}

\theoremstyle{remark}

\renewcommand{\emptyset}{\varnothing}
\newcommand{\CC}{{\mathbb{C}}}
\newcommand{\RR}{{\mathbb{R}}}
\newcommand{\ZZ}{{\mathbb{Z}}}
\newcommand{\QQ}{{\mathbb{Q}}}
\newcommand{\PP}{{\mathbb{P}}}
\newcommand{\NN}{{\mathbb{N}}}
\newcommand{\down}[1]{\left\lfloor #1\right\rfloor}
\newcommand{\ep}{\varepsilon}
\newcommand{\var}{\varphi}
\newcommand{\ov}[1]{\overline{#1}}
\renewcommand{\tilde}[1]{\widetilde{#1}}
\newcommand{\rk}{\operatorname{rk}}
\newcommand{\Weil}{\operatorname{Weil}}
\newcommand{\ver}{\operatorname{ver}}
\newcommand{\hor}{\operatorname{hor}}
\newcommand{\Pic}{\operatorname{Pic}}
\newcommand{\Diff}{\operatorname{Diff}}
\newcommand{\Sing}{\operatorname{Sing}}
\newcommand{\pt}{\operatorname{pt}}
\newcommand{\Supp}{\operatorname{Supp}}
\newcommand{\OOO}{{\mathcal{O}}}
\newcommand{\EEE}{{\mathcal{E}}}
\newcommand{\LLL}{{\mathcal{L}}}
\newcommand{\Ga}{\Gamma}

\begin{document}

\title{On a conjecture of Shokurov: \\
Characterization of toric varieties}
\author{Yuri G. Prokhorov}
\address{Department of Mathematics
(Algebra Section), Moscow Lomonosov University, 117234 Moscow,
Russia\quad \&\quad Department of Mathematics, Tokyo Institute of
Technology, Oh-okayama, Meguro, 152 Tokyo, Japan}

\email{prokhoro@mech.math.msu.su}
\begin{abstract}
We verify a special case of V. V. Shokurov's conjecture about
characterization of toric varieties. More precisely, let
$(X,D=\sum d_iD_i)$ be a three-dimensional log variety such that
$K_X+D$ is numerically trivial and $(X,D)$ has only purely log
terminal singularities. In this situation we prove the inequality
\begin{center}
$\sum d_i\le \rk\Weil(X)/(\operatorname{algebraic\ equivalence})
+\dim(X)$.
\end{center}
We describe such pairs for which the equality holds and show that
all of them are toric.
\end{abstract}
\maketitle
\section{Introduction}
The aim of this note is to discuss the birational characterization
of toric varieties. Let $X$ be a normal projective toric variety
and let $D=\sum_{i=1}^rD_i$ be the sum of invariant divisors. It
is well known that the pair $(X,D)$ has only log canonical
singularities (see e. g. \cite[3.7]{Ko}), $K_X+D$ is linearly
trivial and $r=\rk(\Weil(X)/\approx)+\dim(X)$, where $\Weil(X)$ is
the group of Weil divisors and $\approx$ is the algebraic
equivalence.

Shokurov observed that this property can characterize toric
varieties:

\begin{conjecture}[{\cite{Sh1}}]
\label{conjecture}
Let $(X,D=\sum d_iD_i)$ be a projective log variety such that
$(X,D)$ has only log canonical singularities and numerically
trivial. Then $\sum d_i\le \rk(\Weil(X)/\approx)+\dim(X)$.
Moreover, if the equality holds, then $(X,\down{D})$ is a toric
pair.
\end{conjecture}

Shokurov also conjectured the relative version of \ref{conjecture}
(cf. Theorem~\ref{local}) and he expects that one can replace
numerical triviality of $K_X+D$ with nefness of $-(K_X+D)$. We do
not discuss these details here.

Conjecture~\ref{conjecture} was proved in dimension two in
\cite{Sh1} (see also \cite[Sect. 8]{Lect} and
Proposition~\ref{Conj-2} below). Our main result is
the following partial answer on 
Conjecture~\ref{conjecture} in dimension three:

\begin{theorem}
\label{Th-1}          4
Let $(X,D=\sum d_iD_i)$ be a three-dimensional projective variety
over $\CC$ such that $K_X+D\equiv 0$ and $(X,D)$ has only purely
log terminal singularities. Then
\setcounter{equation}{\value{subsection}}
\begin{equation}
\sum d_i\le \rk(\Weil(X)/\approx)+3.
\label{eq-main-Weil}
\end{equation}
\setcounter{subsection}{\value{equation}} Moreover, if the
equality holds, then up to isomorphisms one of the following
holds:
\begin{enumerate}
\renewcommand\labelenumi{(\roman{enumi})}
\item
$X\simeq\PP^3$, $\down{D}=0$ or $\down{D}=\PP^2$;
\item
$X\simeq\PP^1\times\PP^2$, $\down{D}=0$ or $\down{D}=\{\pt\}\times\PP^2$ or 
$\down{D}=\PP^1\times\{\operatorname{line}\}$;
\item
$X\simeq\PP\left(\OOO_{\PP^1}\oplus \OOO_{\PP^1}\oplus
\OOO_{\PP^1}(d)\right)$, $d\ge 1$, $\down{D}$ is the section correspondong 
to the surjection $\OOO_{\PP^1}\oplus \OOO_{\PP^1}\oplus
\OOO_{\PP^1}(d)\to \OOO_{\PP^1}(d)$;
\item
$X\simeq\PP^1\times\PP^1\times\PP^1$, $\down{D}=0$ 
or $\down{D}=\{\pt\}\times\PP^1\times\PP^1$ 
or $\down{D}=\{\pt_1,\pt_2\}\times\PP^1\times\PP^1$;
\item
$X\simeq\PP\left(\OOO_{\PP^2}\oplus \OOO_{\PP^2}(d)\right)$, $d\ge 1$, 
$\down{D}$ is the negative section, or a disjoint
union of two sections, one of them is negative;
\item
$X\simeq\PP\left(\OOO_{\PP^1\times\PP^1}\oplus
\LLL\right)$, $\LLL\in\operatorname{Pic}(\PP^1\times\PP^1)$, 
$\down{D}$ is the negative section, or a disjoint
union of two sections, one of them is negative.
\end{enumerate}
In all cases $(X,\down{D})$ is toric.
\end{theorem}
Clearly, our theorem is not a characterization of toric varieties,
but we hope that Conjecture~\ref{conjecture} can be proved in a
similar way.

This paper was subject of my talk given at Waseda seminar on
January 18, 2000. I am grateful to the participants of this
seminar, especially Professor S. Ishii,
for attention and valuable discussions. 
I also would like to thank the
Department of Mathematics of Tokyo Institute of Technology for
hospitality during my stay on 1999--2000. This work was partially
supported by the grant INTAS-OPEN-97-2072.

\section{Preliminaries}
\subsection*{Notation}
All varieties are defined over $\CC$. Basically we employ the
standard notation of the Minimal Model Program (MMP, for short).
Throughout this paper $\rho(X)$ is the Picard number and
$\ov{NE}(X)$ is the Mori cone of $X$. We call a pair $(X,D)$
consisting of a normal algebraic variety $X$ and a boundary $D$ on
$X$ a \textit{log variety} or a \textit{log pair}. Here a
\textit{boundary} is a $\QQ$-Weil divisor $D=\sum d_iD_i$ such
that $0\le d_i\le 1$ for all $i$. A \textit{contraction} is a
projective morphism $\var\colon X\to Z$ of normal varieties such
that $\var_*\OOO_X=\OOO_Z$. Abbreviations klt, plt, lc are
reserved for Kawamata log terminal, purely log terminal and log
canonical, respectively (refer to \cite{Sh}, \cite{Ut} and
\cite{Ko} for the definitions). Let $(X,D)$ be a log pair and let
$S:=\down{D}$. For simplicity, assume that $(X,D)$ is lc in
codimension two. The Adjunction Formula proposed by Shokurov
\cite[Sect. 3]{Sh} states $(K_X+D)|_S=K_S+\Diff_S(D-S)$, where
$\Diff_S(D-S)$ is a naturally defined effective $\QQ$-Weil divisor
on $S$, so-called \textit{different}. Moreover, $K_X+D$ is plt
near $S$ if and only if $S$ is normal and $K_S+\Diff_S(D-S)$ is
klt \cite[17.6]{Ut}. $LCS(X,D)$ denotes the \textit{locus of log
canonical singularities} of $(X,D)$ that is the set of all points
where $ (X,D)$ is not klt \cite{Sh}. Let $\var\colon X\to Z$ be
any fiber type contraction and let $D=\sum d_iD_i$ be a
$\QQ$-divisor on $X$. We will write
$D=\sum_{\ver}d_iD_i+\sum_{\hor}d_iD_i= D^{\ver}+D^{\hor}$, where
$\sum_{\ver}$ (respectively $\sum_{\hor}$) runs through all
components $D_i$ such that $\dim \var (D_i)<\dim(D_i)$
(respectively $\var (D_i)=Z$). We will frequently use the above
notation without reference.

In dimension two Conjecture~\ref{conjecture} is much easier than
higher dimensional one. We need only the following weaker version:

\begin{proposition}
\label{Conj-2}Let $(X,D=\sum d_iD_i)$ be a projective log surface such that $
-(K_X+D)$ is nef and $(X,D)$ is lc. Then $\sum d_i\le\rho(X)+2$.
Moreover, if the equality holds and $(X,D)$ is klt, then
$X\simeq\PP^2$, or $X\simeq\PP^1\times\PP^1$.
\end{proposition}
For the general statement we refer to \cite{Sh1}, see also
\cite{Lect}.
\begin{proof} Assume that $\sum d_i-\rho(X)-2\ge 0$ and run $K_X$-MMP.
According to \cite{A}, Log MMP works even in the category of log
canonical pairs. On each step, $\sum d_i-\rho (X)-2$ does not
decrease and all assumptions are preserved (see \cite[2.28]{Ut}).
At the end we get one of the following:

\subsection*{Case~1} $\rho(X)=1$. Then $\sum d_i-\rho(X)-2\le 0$ by
\cite[18.24]{Ut}.

\subsection*{Case~2}
There is an extremal contraction onto a curve $\var\colon X\to Z$
(in particular, $\rho(X)=2$). Let $\ell$ be a general fiber. Then
\setcounter{equation}{\value{subsection}}
\begin{equation}
\label{2-2-2}
2=-K_X\cdot\ell\ge D\cdot\ell\ge\sum_{\hor}d_i.
\end{equation}
\setcounter{subsection}{\value{equation}} Hence $\sum_{\ver}d_i\ge
2$ and $K_X+D^{\hor}$ is not nef. Let $\phi\colon X\to W$ is a
contraction of $(K+D^{\hor})$-negative extremal ray. If $\phi$ is
birational, we replace $X$ with $W$ and obtain Case~1 above. Thus
we may assume that $W$ is a curve, so $\var$ and $\phi$ are
symmetric. As above, $\sum_{\ver}d_i\le 2$, so $\sum d_i=4$.

Now assume that $(X,D)$ is klt and $\sum d_i-\rho(X)-2=0$. Then
after each divisorial contraction $\sum d_i-\rho(X)-2$ increases.
Hence we have only cases 1 or 2 above. In Case~1, $X\simeq\PP^2$
by Lemma~\ref{rho=1} below. In Case~2 we have the equality in
(\ref{2-2-2}), so $D_i\cdot \phi^{-1}(w)=1$ for any component of
$D^{\ver}$ and a general fiber $\phi^{-1}(w)$. Therefore $D_i$ is
not a multiple fiber of $\var$ and $X$ is smooth along $D_i$.
Similarly, if $D_j$ is a component of $D^{\hor}$, then $\phi
(D_j)=\pt$, $X$ is smooth along $D_j$ and $D_i\cdot D_j=1$.
Finally, $\var \times\phi \colon X\to Z\times W=\PP^1\times\PP^1$
is a finite morphism of degree
$\var^{-1}(\pt)\times\phi^{-1}(\pt)=D_i\cdot D_j=1$.
\end{proof}

The local version of Conjecture~\ref{conjecture} was proved in
\cite[18.22]{Ut}:

\begin{theorem}
\label{local}Let $(X,D=\sum d_iD_i)$ be a log pair which is log canonical at
a point $P\in \cap D_i$. Assume that $K_X$ and all $D_i$'s are
$\QQ$-Cartier at $P$. Then $\sum d_i\le \dim(X)$. Moreover, if the
equality holds, then $(X\ni P,D)$ is an abelian quotient of a
smooth point and $(X,D)$ is not plt at $P$.
\end{theorem}

Recall that for any plt pair $(X,D)$ of dimension $\le 3$ there is
a small birational contraction $q\colon X^q\to X$ such that $X^q$
is $\QQ $-factorial and $(X^q,D^q:=q_*^{-1}D)$ is plt (see
\cite[6.11.1]{Ut}, \cite[17.10]{Ut}). Such $q$ is called a
$\QQ$\textit{-factorialization} of $(X,D)$. Applying a
$\QQ$-factorialization in our situation and taking into account
that $\rk\left(\Weil(X)/\approx\right)
=\left(\rk\Weil(X^q)/\approx\right) \ge\rho(X^q)$ we obtain that
for Theorem~\ref{Th-1} it is sufficient to prove the following

\begin{proposition}
\label{Th-2}
Let $(X,D=\sum d_iD_i)$ be a three-dimensional projective plt
pair such that $K_X+D\equiv 0$ and $X$ is $\QQ$-factorial. Then
\setcounter{equation}{\value{subsection}}
\begin{equation}
\sum d_i\le\rho(X)+3.
\label{eq-main-rho}
\end{equation}
\setcounter{subsection}{\value{equation}} Moreover, if the
equality holds, then for $(X,\down{D})$ there are only
possibilities {\rm (i)--(iv)} of Theorem~\ref{Th-1}.
\end{proposition}

\section{Lemmas}
In this section we prove several facts related to
Conjecture~\ref{conjecture}.

\begin{lemma}[cf. {\cite[18.24]{Ut}}, {\cite{A}}]
\label{rho=1}
Let $(X,D=\sum_{i=1}^rd_iD_i)$ be a projective $n$-dimensional log
pair such that all $D_i$'s are $\QQ$-Cartier, $\rho(X)=1$, $(X,D)$
is plt and $-(K_X+D)$ is nef. Then
\setcounter{equation}{\value{subsection}}
\begin{equation}
\sum d_i\le n+1. \label{eq-rho=1}
\end{equation}
\setcounter{subsection}{\value{equation}} Moreover, if the
equality holds, then $X\simeq\PP^n$ and $D_1,\ldots,D_r$ are
hyperplanes.
\end{lemma}

Note that in the two-dimensional case any plt pair is
automatically $\QQ $-factorial.

\begin{proof}We will prove this lemma in the case when $\down{D}=0$
(i.e. $K_X+D$ is klt). The case when $\down{D}$ is non-trivial
(and irreducible) can be treated in a similar way. The inequality
(\ref{eq-rho=1}) was proved in \cite[18.24]{Ut}, so we prove the
second part of our lemma. Since $-K_X$ is ample,
$\Pic(X)\simeq\ZZ$ (see e.g. \cite[2.1.2]{IP}). Let $H$ be an
ample generator of $\Pic(X)$ and let $D_i\equiv a_iH$,
$a_i\in\QQ$, $a_i\ge 0$. Assume that $a_i<a_j$ for $i\neq j$ or
$K_X+D\not\equiv 0$. For $0<\ep \ll 1$, consider
\[
D^{(\ep)}:=\ep D_i+D-\ep D_j.
\]
Then $K_X+D^{(\ep)}$ is again klt (because the klt property is an
open condition) and $-(K_X+D^{(\ep)})$ is ample. Take $N\in\NN$ so
that $-N(K_X+D^{(\ep)})$ is integral and very ample and let $M\in
|-N(K_X+D^{(\ep)})|$ be a general member. By Bertini type theorem
\cite[Sect. 4]{Ko}, $(X,D^{(\ep)}+\frac 1NM)$ is klt (and
numerically trivial). Moreover, the sum of coefficients of
$D^{(\ep)}+\frac 1NM$ is equal to $n+1+\frac 1N$. This contradicts
(\ref{eq-rho=1}). Hence $K_X+D\equiv 0$ and $D_i\equiv D_j$ for
all $i,j$. Thus, for any pair $D_{i,}D_j$ there exists
$n_{i,j}\in\NN$ such that $n_{i,j}(D_i-D_j)\sim 0$. By taking
repeated cyclic covers (which are \'etale in codimension one)
$\pi\colon X'\to \cdots\to X$, we obtain a new plt pair
$(X',D'=\sum_{i=1}^{r}d_iD_i')$ \cite[20.4]{Ut} such that
$D_i'\sim D_j'$, where $D_i'=\pi^*D_i$. On this step, we do not
assume that $D_i'$ is irreducible. Then $D_1',\ldots,D_{r}'$
generate a linear system $\LLL$ of Weil divisors. If
$\operatorname{Bs}(\LLL)$ is not empty, then we pick a point
$P'\in D_1'\cap \cdots \cap D_{r}'$. By construction, $({X'},D')$
is klt at $P'$ and $\sum_{i=1}^{r}d_i\ge n+1$, a contradiction
with Theorem~\ref{local}. Therefore
$\operatorname{Bs}(\LLL)=\emptyset$. In particular,
$D_1',\ldots,D_{r}'$ are ample Cartier divisors and $-K_{X'}\equiv
D'$ is ample (i.e. $X'$ is a log Fano variety). This also shows
that the Fano index of $X'$ is $r(X')\ge \sum_{i=1}^{r}d_i\ge
n+1$. It is well known (see e.g. \cite[3.1.14]{IP}) that in this
case we have $r(X')=\sum_{i=1}^{r}d_i=n+1$, $X'\simeq\PP^n$ and
$D_1',\ldots,D_{r}'$ are hyperplanes. Since $\pi\colon X'\to X$ is
\'etale outside of $\Sing(X)$ and $X'$ is smooth, the restriction
$X'\backslash \pi^{-1}\left(\Sing(X)\right)\to X\backslash
\Sing(X)$ is the universal covering. This gives us that $\pi\colon
X'\to X$ is Galois. Hence $X=\PP^n/G$, where $G\subset PGL_{n+1}$
is a finite subgroup. Further, the group $G$ does not permute
$D_1',\ldots,D_{r}'$. Thus $G$ has $r>n+1$ invariant hyperplanes
$D_1',\ldots,D_{r}'$ in $\PP^n$. Finally, the lemma follows by the
following simple fact which can be proved by induction on $n$.
\end{proof}

\begin{sublemma*}
Let $G\subset PGL_{n+1}$ be a finite subgroup acting on $\PP^n$
free in codimension one. Assume that there are $n+2$ invariant
hyperplanes $H_1,\ldots,H_{n+2}\subset\PP^n$. Then $G=\{1\}$.
\end{sublemma*}

\begin{lemma}
\label{flip}
Let $\var\colon X\to Z\ni o$ be a three-dimensional flipping
contraction and let $D=\sum d_iD_i$ be a boundary on $X$ such that
$(X,D)$ is plt, $\rho(X/Z)=1$, $-(K_X+D)$ is $\var$-nef and all
$D_i$'s are $\var$-ample. Assume that $X$ is $\QQ$-factorial. Then
$\sum d_i<2$.
\end{lemma}

\begin{proof}
Let $\chi\colon X\stackrel{\var}{\longrightarrow}Z
\stackrel{\phantom{{}^{+}}\var^{+}}{\longleftarrow}X^{+}$ be the
flip with respect to $K_X$ and let $D^{+}=\sum d_iD_i^{+}$ be the
proper transform of $D$. Then all $D_i^{+}$'s are anti-ample over
$Z$. Hence $\var^{+-1}(o)$ is contained in $\cap D_i^{+}$.
Consider a general hyperplane section $H\subset X^{+}$. Then
$(H,D|_H)$ is plt \cite[Sect. 4]{Ko}. Applying Theorem~\ref{local}
to $H$ we obtain $\sum d_i<2$.
\end{proof}

\begin{lemma}
\label{surface} Let $\var\colon X\to Z$ be a contraction from a
projective $\QQ$-factorial threefold onto a surface such that
$\rho(X/Z)=1$. Let $D=\sum d_iD_i$ be a boundary on $X$ such that
$(X,D)$ is lc, $(X,D-\down{D})$ is klt and $-(K_X+D)$ is nef.
Assume that $\down{D}$ has a component $S$ which is generically
section of $\var$. Then $\sum_id_i\le\rho(X)+3$. Moreover, if the
equality holds and $(X,D)$ is plt, then $X$ is smooth, $\var$ is a
$\PP^1$-bundle, $\var|_S\colon S\to Z$ is an isomorphism and
$Z\simeq\PP^2$ or $Z\simeq\PP^1\times\PP^1$.
\end{lemma}

\begin{proof}Assume that $\sum_id_i\ge\rho(X)+3$. Since $-K_X$
is $\var$-ample, a general fiber $\ell$ of $\var$ is isomorphic to
$\PP^1$. We have \setcounter{equation}{\value{subsection}}
\begin{equation}
2=-K_X\cdot\ell =D^{\hor}\cdot\ell\ge\sum_{\hor}d_i,\qquad
\sum_{\ver}d_i\ge\rho(X)+1=\rho(Z)+2.
\label{eq-eq}
\end{equation}
\setcounter{subsection}{\value{equation}} Let $\mu :=\var|_S$.
Write $\Diff_S(D-S)= =\sum_i\beta_i\Theta_i$. Then
\setcounter{equation}{\value{subsection}}
\begin{equation}
\beta_i=1-\frac 1{m_i}+\frac 1{m_i}\sum_{j\in
\mathfrak{M}_i}d_jk_{i,j},
\label{eq-m}
\end{equation}
\setcounter{subsection}{\value{equation}} where $m_i\in\NN\cup
\{\infty \}$, $k_{i,j}\in\NN$ and the sum runs through the set
$\mathfrak{M}_i$ of all components $D_j$ containing $\Theta_i$
(see \cite[3.10]{Sh}). Here $m_i=\infty $ when $(X,D)$ is not plt
along $\Theta_i$. It is easy to see that $\beta_i\ge\sum_{j\in
\mathfrak{M}_i}d_j$. Put $\Xi :=\mu_*\Diff_S(D-S)$ and let $\Xi
=\sum \gamma_i\Xi_i$. Since $-(K_S+\Theta)$ is nef, $(Z,\Xi)$ is
lc \cite[2.28]{Ut}. For any component $D_i$ of $D^{\ver}$ we have
at least one component $ \Theta_j\subset D_i\cap S$ such that $\mu
(\Theta_j)\neq \pt$. This yields
\[
\sum_i\gamma_i\ge\sum_{\mu (\Theta_i)\neq \pt}\beta_i\ge
\sum_{\ver}d_j\ge\rho(Z)+2.
\]
Applying Proposition~\ref{Conj-2} to $(Z,\Xi)$, we obtain
equalities \setcounter{equation}{\value{subsection}}
\begin{equation}
\sum_i\gamma_i=\sum_{\mu (\Theta_i)\neq
\pt}\beta_i=\sum_{\ver}d_j=\rho(Z)+2.
\label{eq-=gd}
\end{equation}
\setcounter{subsection}{\value{equation}} Hence $\sum_{\hor}d_i=2$
and $\sum_id_i=\rho(X)+3$. This shows the first part of the lemma.

Now assume that $(X,D)$ is plt. By Adjunction \cite[17.6]{Ut},
$(S,\Diff_S(D-S))$ is klt and so is $(Z,\Xi)$. Again, by
Proposition \ref{Conj-2} we have either $Z\simeq\PP^2$ or
$Z\simeq\PP^1\times\PP^1$. There exists a \textit{standard form}
of $\var$ (see \cite{Sarkisov}), i.~e. the commutative diagram
\[
\begin{array}{ccc}
\tilde{X} & \dasharrow & X \\ \downarrow & & \downarrow \\
\tilde{Z} & \stackrel{\sigma}{\longrightarrow} & Z
\end{array}
\]
where $\sigma\colon \tilde{Z}\to Z$ is a birational morphism of
smooth surfaces, $\tilde{X}\dasharrow X$ is a birational map and
$\tilde{\var}\colon \tilde{X}\to \tilde{Z}$ is a \textit{standard}
conic bundle (in particular, $\tilde{X}$ is smooth and
$\rho\left(\tilde{X}/X\right)=1$). Take the proper transform
$\tilde{S}$ of $S$ on $\tilde{X}$. For a general fiber
$\tilde{\ell}$ of $\tilde{\var}$ we have
$\tilde{S}\cdot\tilde{\ell}=1$. Since
$\rho(\tilde{X}/\tilde{Z})=1$, $\tilde{S}$ is
$\tilde{\var}$-ample. It gives us that each fiber of
$\tilde{\var}$ is reduced and irreducible, i.~e. the morphism
$\tilde{\var}$ is smooth. By \cite {Iskovskikh}, there exists a
standard conic bundle $\widehat{\var}\colon \widehat{X}\to Z$ and
a birational map $\widehat{X} \dasharrow X$ over $Z$. This map
indices an isomorphism $(\widehat{X}/
\widehat{\var}^{-1}(\mathfrak{M})) \simeq
(X/\var^{-1}(\mathfrak{M}))$, where $\mathfrak{M}\subset Z$ is a
finite number of points. Since both $\var$, $\widehat{\var}$ are
projective and $\rho(X/Z)=\rho (\widehat{X}/Z)=1$, we have
$\widehat{X}\simeq X$. But then $\var \colon X\to Z$ is smooth,
i.e. $\var$ is a $\PP^1$-bundle.

Now we claim that $\mu$ is an isomorphism. Indeed, otherwise $S$
contains a fiber, say $\ell_0$. Then $S$ intersects all
irreducible components of $D^{\hor}-S$. If some component $D_k$ of
$D^{\hor}-S$ does not contain $\ell_0$, then $\var (S\cap D_k)$ is
a component of $\Xi$. By (\ref{eq-=gd}) we have
\[
\rho(Z)+2=\sum_i\gamma_i=\sum_{\mu (\Theta_i)\neq \pt}\beta_i\ge
d_k+\sum_{\ver}d_j>\rho(Z)+2,
\]
which is impossible. Therefore all components of $D^{\hor}$
contain $ \ell_0$. Taking a general hyperplane section as in the
proof of Lemma~\ref{flip}, we derive a contradiction.
\end{proof}

\begin{corollary}
$S$ does not intersects $\Supp(D^{\hor}-S)$ and all components of
$D^{\hor}-S$ are sections of $\var$.
\end{corollary}

\begin{lemma}
\label{eq-rho=2-hor}
Let $\var\colon X\to Z$ be a contraction from a $\QQ$-factorial
threefold onto a curve and let $D=\sum d_iD_i$ be a boundary on
$X$ such that $(X,D)$ is lc, $(X,D-\down{D})$ is klt. Let $F$ be a
general fiber. Assume that $-(K_X+D)$ is $\var$-nef and
$\rho(X/Z)=1$. Then $\sum_{\hor}d_i\le 3$. Moreover, if the
equality holds and $(X,D)$ is plt, then $F\simeq\PP^2$ and for any
component $D_i$ of $D^{\hor}$ the scheme-theoretic restriction
$D_i|_F$ is a line.
\end{lemma}

\begin{proof}
Put $\Delta :=D|_F$. Then $(F,\Delta)$ is lc, $(F,\Delta
-\down{\Delta})$ is klt (see \cite[Sect. 4]{Ko}) and
$-(K_F+\Delta)$ is nef. Moreover, if $(X,D)$ is plt, then so is
$(F,\Delta)$. Write $\Delta =\sum\delta_i\Delta_i$, where all
$\Delta_i$'s are irreducible curves on $F$. Clearly
$D^{\ver}|_F=0$ and $\sum\delta_i\ge\sum_{\hor}d_i$. If $\rho
(F)=1$, then the assertion of \ref{eq-rho=2-hor} follows by
Proposition~\ref{Conj-2}. Assume that $\rho(F)>1$. Let $C$ be an
extremal $K_F$-negative curve on $F$ (note that $K_F$ is not nef).
Then $C$ intersects all components of $\Delta$ (because $\rho
(X/Z)=1$). Let $\upsilon\colon F\to F'$ be the contraction of $C$.
If $F'$ is a curve, then we take $C$ to be a general fiber of
$\upsilon$. By Adjunction, $2=-\deg K_C\ge \deg\Delta|_C$. This
gives us $2\ge\sum\delta_i\ge\sum_{\hor}d_i$. If $ \upsilon$ is
birational, then $(F',\upsilon (\Delta))$ is lc and all components
of $\upsilon (\Delta)$ pass through the point $\upsilon (C)$. By
Theorem~\ref{local}, the sum of coefficients of $\upsilon
(\Delta)$ is $\le 2$. Hence $\sum_{\hor}d_i\le\sum \delta_i\le 3$.
If $(F,\Delta)$ is klt, then so is $(F',\upsilon (\Delta))$ and
the inequality above is strict. Finally, if $(F,\Delta)$ is plt
and $\down\Delta\neq 0$, then we take $C$ to be $(K_F+\Delta
-\down\Delta)$-negative extremal curve. Then $C$ is not a
component of $\down\Delta$. By \cite[3.10]{Ko}, $(F',\upsilon
(\Delta))$ is plt. Again, by Theorem~\ref{local} the sum of
coefficients of $\upsilon(\Delta)$ is strictly less than $2$. So,
$\sum_{\hor}d_i\le\sum\delta_i<3$. This proves Lemma
\ref{eq-rho=2-hor}.
\end{proof}

\begin{corollary}
\label{corollary}
Notation as in Lemma~\ref{eq-rho=2-hor}. Assume additionally that
$X$ is projective, $-(K_X+D)$ is nef (not only over $Z$),
$\sum_{\hor}d_i=3$, $\sum_{\ver}d_i=2$ and $(X,D)$ is plt. If
$\down{D^{\hor}}\neq \emptyset$, then $\down{D^{\hor}}=\down{D}
\simeq\PP^1\times\PP^1$ and $X$ is smooth along $\down{D}$. In
particular, $X$ has at most isolated singularities.
\end{corollary}

\begin{proof}
Put $S:=\down{D}$. By \cite[17.5]{Ut}, $S$ is normal. Since
$\rho(X/Z)=1$, $S$ is irreducible and all components of $D-S$ meet
$S$. Let $\Diff_S(D-S)=\sum \beta_i\Theta_i$. Clearly,
$-(K_S+\Diff_S(D-S))$ is nef. By \cite[17.6]{Ut}, $(S,\Theta)$ is
klt. As in the proof of Lemma \ref{surface} we see $\sum
\beta_i\ge\sum d_i-1\ge 4$. If $\rho(S)=2$, then equalities $\sum
\beta_i=\sum d_i-1=4$ and Proposition~\ref{Conj-2} give us the
assertion. Assume that $\rho(S)>2$. Then some fiber of
$\var|_S\colon S\to Z$ is not irreducible. Let $\Ga$ be its
irreducible component and let $\upsilon\colon S\to S'$ be the
contraction of $\Ga$. Taking into account that $\Ga $ intersects
all components of $D^{\hor}$. As in Lemma~\ref{surface} we get a
contradiction.
\end{proof}

\begin{lemma}[cf. {\cite[6.9]{Sh}}]
\label{lemma-connect}
Let $\var\colon X\to Z\ni o$ be a $K_X$-negative contraction from
a $\QQ$-factorial variety $X$ such that $\rho(X/Z)=1$ and every
fiber has dimension one. Let $D$ be a boundary on $X$ such that
$(X,D-\down{D})$ is klt and $K_X+D$ is $\var$-numerically trivial.
Assume that $\down{D}$ is disconnected near $\var^{-1}(o)$. Then
$K_X+D$ is plt near $\var^{-1}(o)$.
\end{lemma}

\begin{proof}
Regard $\var\colon X\to Z\ni o$ as a germ near $\var^{-1}(o)$. Put
$S:=\down{D}$. Clearly, for a general fiber $\ell$ of $\var$ we
have $-K_X\cdot\ell =D\cdot\ell =2$. If $S'$ is an irreducible
component of $S$ such that $S'\cdot \ell =0$, then
$S'=\var^{-1}(C)$ for a curve $C\subset Z$. In this case, $S'$
contains $\var^{-1}(o)$ and $S$ is connected near $\var^{-1}(o)$.
Therefore $S$ has exactly two connected components $S_1$, $S_2$,
they are irreducible and $S_1\cdot\ell =S_2\cdot\ell =1$. Then
$S_i$, $i=1,2$ meets all components of $\var^{-1}(o)$. Hence
$S_i\cap\var^{-1}(o)$ is $0$-dimensional. Since $Z$ is normal and
$\var|_{S_i}\colon S_i\to Z$ is birational, $S_i\simeq Z$ and
$S_i\cap\var^{-1}(o)$ is a single point. In particular,
$\var^{-1}(o)$ is irreducible. Clearly, $LCS(X,D)\subset S=S_1\cup
S_2$. Assume that $(X,D)$ is not plt. Then there is a divisor $E$
of the function field $K(X)$ with discrepancy $a(E,D)<-1$. Let
$V\subset X$ be its center. Then $V\subset S$ and we may assume
that $V\subset S_1$ (and $V\ne S_1$). Let $L\subset Z$ be any
effective prime divisor containing $\var(V)$ and let
$F:=\var^{-1}(L)$. Clearly, $(X,D+F)$ is not lc near $V$. For
sufficiently small positive $\ep$ the log pair $(X,D+F-\ep S_1)$
is not lc near $V$ and not klt near $S_2$. This contradicts
Connectedness Lemma~\cite[17.4]{Ut}
\end{proof}

\section{Proof of Theorem \ref{Th-1}}
In this section we prove Proposition~\ref{Th-2}.

\subsection{Inductive hypothesis}
Notation and assumption as in Proposition~\ref{Th-2}. Our proof is
by induction on $\rho(X)$. In the case $\rho(X)=1$, the assertion
is a consequence of Lemma~\ref{rho=1}. To prove \ref{Th-2} for
$\rho (X)\ge 2$ we fix $\rho\in\NN$, $\rho >1$. Assume that
inequality (\ref{eq-main-rho}) holds if $\rho(X)<\rho$ and for
$\rho(X)=\rho$ we have \setcounter{equation}{\value{subsection}}
\begin{equation}
\sum d_i-\rho(X)-3\ge 0.
\label{eq-rho}
\end{equation}
\setcounter{subsection}{\value{equation}}
\subsection{}
If $(X,D)$ is klt, then we run $K_X$-MMP. On each step $K\equiv
-D$ cannot be nef. Obviously, all steps preserve our assumptions
(see \cite[2.28]{Ut}) and the left hand side of (\ref{eq-rho})
does not decrease. Moreover, by our assumptions we have no
divisorial contractions on $X$ (because after any divisorial
contraction the left hand side of (\ref{eq-rho}) decreases).
Therefore after a number of flips, we obtain a fiber type
contraction $\var\colon X\to Z$. Since $\rho(X)=\rho\ge 2$,
$\dim(Z)=1$ or $2$. Note that all varieties from
Theorem~\ref{Th-1} have no small contractions. Thus, it is
sufficient to prove Proposition~ \ref{Th-2} on our new model
$(X,D)$.

This procedure does not work if $(X,D)$ is not klt. The difference
is that contractions of components of $D$ do not contradict the
inductive hypothesis. If $(X,D)$ is not klt, then we run
$(K_X+D-\down{D})$-MMP. Note that $\down{D}$ is normal and
irreducible \cite[17.5]{Ut}. For every extremal ray $R$ we have
$\down{D} \cdot R>0$, so we cannot contract an irreducible
component of $\down{D}$. Therefore after every divisorial
contraction $\sum d_i-\rho(X)$ decrease, a contradiction with our
assumption. Thus, all steps of the MMP are flips. By
\cite[2.28]{Ut}, they preserve the plt property of $K+D$. At the
end we get a fiber type contraction $\var\colon X\to Z$, where
$\dim(Z)<3$ and $\down{D}$ is $\var$-ample (i.e.
$\down{D^{\hor}}\neq 0$). Since $\down{D^{\hor}}$ has a component
which intersects all components of $D^{\ver}$,
$\down{D^{\ver}}=0$.

\subsection{Case: $\dim(Z)=1$}
Then $\rho(X)=2$. By Lemma~\ref{eq-rho=2-hor} and our assumption
(\ref{eq-rho}), we have $\sum_{\hor}d_i\le 3$ and
$\sum_{\ver}d_i\ge 2$. In particular, $D^{\ver}\neq 0$. Components
of $D^{\ver}$ are fibers of $\var$, so they are numerically
proportional. Clearly, the log divisor
$K_X+D^{\hor}\equiv-D^{\ver}$ is not nef and curves in fibers of
$\var$ are trivial with respect to it. Let $Q$ be the extremal
$(K_X+D^{\hor})$-negative ray of $\ov{NE}(X) \subset\RR^2$ and let
$\phi\colon X\to W$ be its contraction. It follows by
Lemma~\ref{flip} that $\phi$ cannot be a flipping contraction. Let
$\ell$ be a general curve such that $\phi(\ell)=\pt$. Then $\ell$
dominates $Z$ and $\ell\simeq\PP^1$. Hence $Z\simeq\PP^1$.

\subsubsection{Subcase $\down{D}=0$}
\label{klt-case=1}
We will prove that $X\simeq\PP^2\times\PP^1$. By our inductive
hypothesis, $\phi$ cannot be divisorial. Therefore $\dim(W)=2$.
Further, $D^{\ver}\cdot \ell\le D\cdot \ell=-K_X\cdot \ell=2$.
Since $\ell$ intersects all components of $D^{\ver}$,
$\sum_{\ver}d_i\le 2$. This yields $\sum_{\ver}d_i=2$ and
$\sum_{\hor}d_i=3$. In particular, this proves inequality
(\ref{eq-main-rho}). Moreover, $\ell \cdot D^{\ver}=2$, $\ell
\cdot D^{\hor}=0$ and for any component $D_i$ of $D^{\ver}$ we
have $D_i\cdot\ell =1$. Fix two components of $D^{\ver}$, say
$D_0$ and $D_1$. Then $K_X+D_0+D_1+D^{\hor}\equiv K_X+D\equiv 0$,
so $(X,D_0+D_1+D^{\hor})$ is plt by Lemma~\ref{lemma-connect}.
Applying Lemma~\ref{surface} we obtain $D_0\simeq D_1\simeq
Z\simeq\PP^2$, $X$ is smooth and $\phi$ is a $\PP^1$-bundle. By
\cite[3.5]{Mo}, $\var$ is a $\PP^2$-bundle. We have a finite
morphism $\var \times\phi\colon X\to Z\times W=\PP^1\times\PP^2$.
Clearly, $\deg (\var \times\phi)=\var^{-1}( \pt)\cdot\ell=1$.
Hence $\var \times\phi$ is an isomorphism.

\subsubsection{Subcase: $\down{D}\ne 0$}
Since $\rho(X/Z)=1$, $\down{D}$ is irreducible. Put $S:=\down{D}$.
Let $F$ be a general fiber of $\var$. By construction, $-K_X$ is
$\var$-ample. First, assume that $\dim(W)=2$. Then
$D^{\ver}\cdot\ell\le D\cdot\ell =-K_X\cdot\ell =2$. Since $\ell$
intersects all components of $D^{\ver}$, $\sum_{\ver}d_i\le 2$.
This yields $\sum_{\ver}d_i=2$, $\sum_{\hor}d_i=3$ and $\sum
d_i=5$. Moreover, $D^{\hor}\cdot\ell =0$. By
Lemma~\ref{eq-rho=2-hor}, $F\simeq\PP^2$, $X$ is smooth along $F$
and for any component $D_i$ of $D^{\hor}$ the scheme-theoretic
restriction $D_i|_F$ is a line. Hence components of $D^{\hor}$ are
numerically equivalent. Let $D_1$ be a component of $D^{\hor}-S$.
Consider the new boundary $D':=D+\ep D_0-\ep S$. If $0<\ep \ll 1$,
then $(X,D')$ is klt and $K_X+D'\equiv 0$. Applying Case
\ref{klt-case=1} we get $W\simeq\PP ^2$ and
$X\simeq\PP^2\times\PP^1$.

Now assume that $\phi$ is divisorial. By the inductive hypothesis,
$\phi$ contract $S$. Since the contraction is extremal, $\phi(S)$
is a curve (otherwise curves $S\cap\var^{-1}(\pt)$ is contracted
by $\phi$ and $\var$). All components of $\phi(D^{\ver})$ pass
through $\phi(S)$. By taking a general hyperplane section as in
the proof of Lemma \ref{flip}, we obtain $\sum_{\ver}d_i\le 2$. By
Corollary \ref{corollary}, we obtain that
$S\simeq\PP^1\times\PP^1$, $X$ has only isolated singularities and
$X$ is smooth along $S$. By Lemma~\ref{eq-rho=2-hor},
$F\simeq\PP^2$ and $X$ is smooth along $F$. The curve $F\cap S$ is
ample on $F$, so it is connected and smooth by the Bertini
theorem. Therefore $F\cap S$ is a is a generator of
$S=\PP^1\times\PP^1$. Since $\var|_S$ is flat, he same holds for
arbitrary fiber $F_0$. Hence all fibers of $\var$ are numerically
equivalent and any fiber $F_0$ contains an ample smooth rational
curve. Moreover, this also means $F_0$ is not multiple. Thus it is
a normal surface. Now as in Case \ref{klt-case=1},
$K_X+F_0+F_1+D^{\hor}\equiv 0$ and by Lemma~\ref{lemma-connect},
$(X,F_0+F_1+D^{\hor})$ is plt for any fibers $F_0$, $F_1$. By
Adjunction, $(F_0,D^{\hor}|_{F_0})$ is klt. Clearly,
$K_{F_0}\equiv K_X|_{F_0}$ and $S|_{F_0}$ are numerically
proportional. Hence $F_0$ is a log del Pezzo surface of Fano index
$>1$. Since $\var$ is flat, $\left(K_{F_0}\right)^2=\left(
K_F\right)^2=9$. Therefore, $F_0\simeq\PP^2$ and $X$ is smooth. By
\cite[3.5]{Mo}, $\var$ is a $\PP^2$-bundle so,
$X\simeq\PP_{\PP^1}(\EEE)$, where $\EEE=\OOO_{\PP^1}\oplus
\OOO_{\PP^1}(a) \oplus \OOO_{\PP^1}(b)$, $0\le a\le b$. The
Grothendiek tautological bundle $\OOO_{\PP(\EEE)}(1)$ is
generated by global sections and not ample. Therefore $\OOO_{\PP(\EEE)}(1)$ gives
us a supporting function for the extremal ray $Q$. Since $\phi$ is
birational, $\OOO_{\PP(\EEE)}(1)^3=a+b>0$. Finally, $X$ contains
$S=\PP^1\times\PP^1$. Hence $a=0$. This proves
Proposition~\ref{Th-2} in the case when $Z$ is a curve.

\subsection{Case: $\dim(Z)=2$}
\label{case-2}
Note that $Z$ has only log terminal singularities (see e. g.
\cite[15.11]{Ut}). Since $-K_X$ is $\var$-ample, a general fiber
$\ell$ of $\var$ is $\PP^1$. Hence $2=-K_X\cdot\ell =D\cdot\ell
=D^{\hor}\cdot\ell \ge\sum_{\hor}d_i$. By our assumption,
$\sum_{\ver}d_i\ge\rho (X)+1$. If $(X,D)$ is not plt, then
$\down{D}$ is $\var$-ample. Clearly, $\down{D^{\ver}}=0$.

\begin{claim}
\label{claim-0}
Notation as above. Then $K_Z$ is not nef.
\end{claim}
\begin{proof}
Run $(K_X+D^{\ver})$-MMP. After a number of flips we get either a
divisorial contraction (of the proper transform of a component of
$\down{D}$), or a fiber type contraction. In both cases $Z$ is
dominated by a family of rational curves \cite[5-1-4, 5-1-8]{KMM}.
Therefore $K_Z$ is not nef by \cite{MM}.
\end{proof}

\begin{claim}
\label{claim}
Notation as in \ref{case-2}. Then $Z$ contains no contractible
curves. In particular, $\rho(Z)\le 2$.
\end{claim}

\begin{proof}
Assume the converse. i. e. there is an irreducible curve
$\Ga\subset Z$ and a birational contraction $\mu\colon Z\to Z''$
such that $\mu(\Ga)=\pt$ and $\rho(Z/Z'')=1$. Denote
$F:=\var^{-1}(\Ga)$. Since $F\not\subset \down{D}$, $(X,D+\ep F)$
is plt for $0<\ep\ll 1$ \cite[2.17]{Ut}. First we assume that
$\Ga$ is contractible: Run $(K+D+\ep F)$-MMP over $Z''$. By our
inductive hypothesis, there are no divisorial contractions
(because such a contraction must contract $F$). At the end we
cannot get a fiber type contraction (because $K+D+\ep F\equiv \ep
F$ cannot be anti-ample over a lower-dimensional variety). Thus
after a number of flips $X\dasharrow X'$, we get a model $X'$ over
$Z''$ such that $K_{X'}+D'+\ep F'\equiv\ep F'$ is nef over $Z''$,
where $D'$ and $F'$ are proper transforms of $D$ and $F$,
respectively. Then $F'\not\equiv 0$ (because $F\not\equiv 0$). Let
$\ell'$ be the proper transform of a general fiber of $\var$.
Since $F'$ is nef over $Z''$, $F'\cdot \ell'=0$ and
$\rho(X'/Z'')=2$, we obtain that $\ell'$ generates an extremal ray
of $\ov{NE}(X'/Z'')$. Let $\var'\colon X'\to Z'$ be its
contraction and let $\mu'\colon Z'\to Z''$ be the natural map.
Then $\dim(Z')=2$, $\Ga':=\var'(F')$ is a curve and
$\mu'(\Ga')=\mu(\var(F))=\mu(\Ga)=\pt$. Therefore
$\left(\Ga'\right)^2<0$. On the other hand,
$\left(\Ga'\right)^2\ge 0$, which is a contradiction. Indeed, Let
$C'\subset F'$ be any curve such that $\var'(C')=\Ga'$. Then
$C'\cdot F'\ge 0$. By the projection formula,
$\left(\Ga'\right)^2\ge 0$.

Since $K_Z$ is not nef, there is an extremal contraction
$\psi\colon Z\to V$. By the above it is not birational. Therefore
$\dim(V)=1$ and $\rho(Z)=2$.
\end{proof}

\begin{corollary}
\label{claim-corollary}
Notation as in \ref{case-2}. One of the following holds:
\begin{enumerate}
\renewcommand\labelenumi{(\roman{enumi})}
\item
$\rho(Z)=1$ and $-K_Z$ is ample;
\item
$\rho(Z)=2$ and there is a $K_Z$-negative extremal contraction
$\psi\colon Z\to V$ onto a curve.
\end{enumerate}
\end{corollary}

\subsubsection{Subcase $\down{D}=0$}
Let $D_i$ be a component of $D^{\ver}$. Run $(K+D-d_iD_i)$-MMP:
\[
\chi_{(i)}\colon X\dasharrow X^{(i)}.
\]
As above we get a fiber type contraction $\var_{(i)}\colon
X^{(i)}\to Z^{(i)}$. Notations $D^{\ver}$ and $D^{\hor}$ will be
fixed with respect to our original $\var$. If $\dim(Z^{(i)})=1$,
then replacing $X$ with $X^{(i)}$ we get the case $\dim(Z)=1$
above. Thus we can assume that $\dim(Z^{(i)})=2$ for any chose of
$D_i$. Let $\ell^{(i)}\subset X^{(i)}$ be a general fiber of
$\var_{(i)}$ and let $L^{(i)}\subset X$ be its proper transform.
Clearly, $\chi_{(i)}$ is an isomorphism along $L^{(i)}$. Hence
$-K_X\cdot L^{(i)}=2$, $L^{(i)}$ is nef and $D_i\cdot L^{(i)}>0$.
For $i=1,\ldots,r$ we get rational curves
$L^{(1)},\ldots,L^{(r)}$. We shift indexing so that $X=X^{(0)}$
and put $Z=Z^{(0)}$ and $\var =\var_{(0)}$.

Up to permutations we can take $L^{(0)},\ldots,L^{(s)}$, $s+1\le
r$ to be linearly independent in $N_1(X)$. Then for any $D_i$
there exists $L^{(j)}$ such that $D_i\cdot L^{(j)}>0$. Thus we
have
\begin{multline*}
2(s+1)=-K_X\cdot\sum_{j=0}^sL^{(j)}=D\cdot\sum_{j=0}^sL^{(j)}=\\
\sum_{i=1}^rd_i\left(D_i\cdot\sum_{j=0}^sL^{(j)}\right)\ge
\sum_{i=1}^rd_i\ge\rho(X)+3.
\end{multline*}
Since $3\ge\rho(Z)+1=\rho(X)\ge s+1$, this yields 
$\rho(X)=s+1=3$. Thus, 
\setcounter{equation}{\value{subsection}}
\begin{equation}
D_i\cdot \sum_{j=0}^2L^{(j)}=1
\label{eq-DL-DL}
\end{equation}
\setcounter{subsection}{\value{equation}} holds for all $i$.
Moreover, $L^{(0)},L^{(1)},L^{(2)}$ generate $N_1(X)$ and
components of $D$ generate $N^1(X)$.

Taking into account that $2=-K_X\cdot L^{(j)}=D\cdot L^{(j)}$, we
decompose $D$ into the sum $D=D^{(0)}+D^{(1)}+D^{(2)}$ of
effective $\QQ$-divisors without common components so that
\setcounter{equation}{\value{subsection}}
\begin{equation}
D^{(i)}\cdot L^{(j)}=\left\{
\begin{array}{lll}
0 & \text{if} & i\neq j \\ 2 & \text{otherwise} &
\end{array}
\right.
\label{eq-DL}
\end{equation}
\setcounter{subsection}{\value{equation}} Then
$D^{(i)}=\var^*\Delta^{(i)}$ for $i=1,2$, where $\Delta^{(1)}$,
$\Delta^{(2)}$ are effective $\QQ$-divisors on $Z$. Put
$C^{(i)}:=\var (L^{(i)})$, $i=1,2$. Since families $L^{(j)}$ are
dense on $X$, $C^{(j)}$ are nef and $\not\equiv 0$. By the
projection formula,
\[
\Delta^{(i)}\cdot C^{(j)}\quad\left\{
\begin{array}{lll}
=0 & \text{if} & 1\le i\neq j\le 2 \\
>0 & \text{if} & 1\le i=j\le 2
\end{array}
\right.
\]
Hence $\Delta^{(1)}$ and $\Delta^{(2)}$ generate extremal rays of
$ \ov{NE}(Z)\subset\RR^2$. By (\ref{eq-DL}), these $\QQ$-divisors
have more than one component, so they are nef and $\left(
\Delta^{(1)}\right)^2=\left(\Delta^{(2)}\right)^2=0$. This gives
us that $C^{(1)}$ and $C^{(2)}$ also generate extremal rays.
Therefore $C^{(i)}$ and $\Delta^{(j)}$ are numerically
proportional whenever $i\neq j$ and $\left(C^{(1)}\right)^2=\left(
C^{(2)}\right)^2=0$. In particular, $C^{(i)}$, $i=1,2$ generate an
one-dimensional base point free linear system which defines a
contraction $Z\to\PP^1$. This also shows that
$D^{(i)}=\var^*\Delta^{(i)}$, $i=1,2$ are nef on $X$.

Now we claim that $D^{(0)}$ is nef. Assume the opposite. Then for
small $\ep>0$, $(X,D+\ep D^{(0)})$ is klt \cite[2.17]{Ut}. There
is a $(K_X+D+\ep D^{(0)})$-negative extremal ray, say $R$. By out
inductive hypothesis, the contraction of $R$ must be of flipping
type. Since $\Delta^{(1)}$, $\Delta^{(2)}$ generate $N^1(Z)$, we
have $D^{(i)}\cdot R>0$ for $i=1$ or $2$. By (\ref{eq-DL}),
$\sum_j^{(i)}d_j=2$, where $\sum_j^{(i)}$ runs through all
components $D_j$ of $D^{(i)}$. Since components of
$D^{(i)}=\var^*(\Delta^{(i)})$ are numerically proportional, this
contradicts Lemma \ref{flip}. Therefore $D^{(i)}$ are nef for
$i=0,1,2$.

We claim that $L^{(i)}$, $i=0,1,2$ generate $\ov{NE}(X)$. Indeed,
let $z\in \ov{NE}(X)$ be any element. Then $z\equiv\sum
\alpha_i[L^{(i)}]$ for $\alpha_i\in\RR$. By (\ref{eq-DL}), $0\le
D^{(j)}\cdot z=\alpha_j$. This shows that $L^{(i)}$ generate
$\ov{NE}(X)$. From (\ref{eq-DL-DL}) we see that components of
$D^{(0)}$ are numerically equivalent.

Fix two components $D'$ and $D''$ of $D^{(0)}$. Then
$K_X+D'+D''+D^{(1)}+D^{(2)}\equiv 0$. By
Lemma~\ref{lemma-connect}, $(X,D'+D''+D^{(1)}+D^{(2)})$ is plt.
Lemma~\ref{surface} implies that $D'\simeq D''\simeq
S\simeq\PP^1\times \PP^1$, $X$ is smooth and $\var$ is a
$\PP^1$-bundle. As in the case $\dim(Z)=1$, we have
$X\simeq\PP^1\times\PP^1\times\PP^1$.

\subsubsection{Subcase $\down{D}\ne 0$}
Let $S$ be a component of $\down{D}$. Clearly, $S\cdot\ell\le 2$.
If $S$ is generically a section of $\var$, then by
Lemma~\ref{surface}, $X$ is smooth, $\var$ is a $\PP^1$-bundle and
$S\simeq\PP^2$ or $\PP^1\times\PP^1$. Therefore $X\simeq
\PP(\EEE)$, where $\EEE$ is a rank two vector bundle on $Z$. Since
$\var$ has disjoint sections, $\EEE$ is decomposable.
So we may assume that $\EEE=\OOO_Z+\LLL$, where $\LLL$ is a line bundle.
By the projection formula, all components of $D^{\ver}$ are nef.
Let $R$ be a $(K_X+D^{\hor})$-negative extremal curve and let 
$\phi\colon X\to W$ be its contraction.
Assume that $\phi$ is of flipping type. 
By \cite{Mo}, $K_X\cdot R\ge 0$.   Hence $D^{\hor}\cdot R<0$, so 
$R$ is contained in a section of $\var$. But all curves on 
$\PP^2$ and $\PP^1\times\PP^1$ are movable, a contradiction.
If $\phi$ is of fiber type, then as in the case $\down{D}=0$
we get $X\simeq Z\times\PP^1$. Assume that $\phi$ is of divisorial type.
By inductive hypothesis, $\phi$ contracts a component of $\down{D}$.

Finally, consider the case when $\var|_S\colon S\to Z$ is generically finite
of degree $2$. Obviously, $D^{\hor}=S$. If $\rho(Z)=1$, then all
components of $D^{\ver}$ are numerically proportional and
$\sum_{\ver}d_i\ge 4$. If $\dim(W)=2$, then
$D^{\ver}\cdot\phi^{-1}(w)\le 2$ for general $w\in W$. Hence
$\sum_{\ver}d_i\ge 4$, a contradiction. Then by lemmas \ref{flip}
and \ref{eq-rho=2-hor}, $\phi$ is divisorial and $\phi$ must
contract $S$. We derive a contradiction with Theorem \ref{local}
for $(W,\phi (D^{\ver}))$ near $\phi (S)$.

Therefore $\rho(Z)=2$ and there is a $K_Z$-negative extremal
contraction $\psi\colon Z\to V$ onto a curve (see
Corollary\ref{claim-corollary}). Let $\pi\colon
X\stackrel{\var}{\longrightarrow}Z
\stackrel{\psi}{\longrightarrow}V$ be the composition map.
Clearly, all fibers of $\pi$ are irreducible. Write
$D=\sum'd_iD_i+\sum''d_iD_i= D'+D''$, where $\sum'$ (respectively
$\sum''$) runs through all components $D_i$ such that
$\pi(D_i)=\pt$ (respectively $\phi(D_i)=V$). Thus, $S$ is a
component of $D''$ and components of $D'$ are numerical
proportional. Let $F$ be a general fiber. Then $\rho(F)=2$.
Consider the contraction $\var|_F\colon F\to \var(F)$ and denote
$D''|_F=D|_F$ by $\Phi=\sum\alpha_i\Phi_i$. Then $(F,\Phi)$ is plt
and $K_F+\Phi\equiv 0$. Clearly, the curve $S|_F=\down{\Phi}$ is a
$2$-section and components of $\Phi-\down{\Phi}$ are fibers of
$\var|_F$. As in the proof of Lemma~\ref{eq-rho=2-hor} using the
fact that $S|_F$ intersects components of $\Phi-\down{\Phi}$
twice, one can check $\sum \alpha_i<3$. This yields $\sum''d_i<3$
and $\sum'd_i>2$. Let $R$ be a $(K_X+D'')$-negative extremal ray.
Since $\sum'd_i>2$ and $\rho(X)>2$, $R$ cannot be fiber type.
According to Lemma~\ref{flip} $R$ also cannot be flipping type.
Therefore $R$ is divisorial and contract $S$ to a point. Since $S$
intersects all components of $D^{\ver}$, this contradicts
Theorem~{local}. The proof of Proposition~\ref{Th-2} is finished.

\subsubsection*{Concluding remark}
In contrast with the purely log terminal case we have no complete
results in the log canonical case. The reason is that the steps of
MMP are not so simple. In particular, we can have divisorial
contractions which contract components of $\down{D}$.

\end{document}